# Non-reflective traveling waves in finite thin beams: A parametric study

Amirhossein Omidi Soroor.[1] and Pablo A. Tarazaga

J. Mike Walker '66 Department of Mechanical Engineering, Texas A&M University

College Station, TX 77843, United States

## Abstract

In the present study, the authors introduce a geometrically improved model inspired by the mammalian basilar membrane's properties and special vibratory behavior while conducting a parametric investigation. The goal of this model is to mimic the broadband non-reflective traveling wave response to excitation frequencies, as observed in the basilar membrane. Simple structural elements, such as beams, springs, and dashpots are utilized for this purpose. Therefore, the beam's equation of motion is developed using Hamilton's principle, and the Galerkin method is implemented as the discretization scheme. Then, the model is verified by comparing its outcome to results presented in the literature. Finally, a parametric study is conducted to reveal the effect of different parameters, i.e., the absorber's location and coefficients' values, excitation frequency, geometric tapering, and material grading on the non-reflective traveling wave response of the beam. The results reveal that traveling waves and their quality are strongly dependent on these parameters. In addition, this study suggests that adopting a single spring-damper system within the span of the beam as the wave reflection absorber may not address the entire bandwidth.

---

[1]Corresponding author.
  *E-mail address:* omidisoroor@tamu.edu (A. Omidi Soroor).

# 1 Introduction

It is well known that when an incident wave propagates through a finite continuum and encounters a boundary or a discontinuity along the way, the wave will reflect, in part or entirely. This, in turn, gives rise to the formation of standing waves (SWs). As most structures are considered finite in real-life instances, SWs dominate the vibratory behavior. However, there are certain cases where traveling waves (TWs) have been observed in natural phenomena. These include but are not limited to certain types of microorganisms such as spermatozoa motility [1], undulatory swimming motion in fishes [2], snail locomotion mechanism [3, 4], and the basilar membrane (BM) dynamic behavior in the mammalian cochlea [5]. Regarding the mammalian ear, it is observed that as the sound signal enters the ear canal, after a series of impedance-matching processes, it excites the BM in the cochlea. Therefore, a wave is initiated in its basal end and propagates toward the apex without any reflection [5]. As the BM stiffness is graded along its length, a TW peaks in magnitude at a location corresponding to the excitation frequency. The higher the excitation frequency, the closer the peak is located to the BM's base [6]. Inspired by the aforementioned observations, there are applications developed where the focus is to leverage this kind of wave for a variety of purposes. Some examples of such applications are particle transfer [7, 8], robot locomotion [9-11], drag reduction [12, 13], and structural health monitoring (e.g., pitch-catch method) [14]. Removing reflections could also improve recent attempts in building occupant tracking [15-17] and classification [18-20] using underfloor accelerometers.

Two major approaches in the literature are concerned with TW generation in finite structures, i.e., two-point and one-point excitation methods. Most published works have leveraged the two-point excitation method to this end. Surveying the literature suggests that Tomikawa et al. [21] were the first to suggest the two-point (or the two-mode) excitation method for TWs generation in 1990. They aimed to develop a more convenient solution to excite TWs for parts-feeder applications than the rather complicated impedance-matching approach discussed in [22]. They suggested that concurrently exciting a finite beam at two different points with a frequency between two adjacent natural frequencies and a 90-degree phase difference could lead to the development of high-quality TWs. This was based on the notion that exciting two consecutive modes of the structure could effectively produce TWs as the result of the interaction of these modes. Since then, there has been a growing interest in the use of the two-point excitation method to induce TWs and its potential applications. For instance, Avirovik et al. [23] studied the two-point excitation model, concluding that any given set of adjacent modes could be used with a 90-degree phase difference to obtain a pure TW response on a beam. Malladi et al.'s [24] study revealed that the presence of TW and its purity is highly dependent on the phase difference. Their study also revealed that based on the excitation frequency, TWs could also be obtained in phase differences other than 90 degrees. Different methods for assessing the quality of TWs were discussed in detail in [25-28].

Of all the natural phenomena discussed above, the BM's dynamic behavior is of particular interest in this work, as a recent study has revealed that unsafe listening practices have endangered more than a billion young individuals' hearing [29]. Therefore, providing more insight into the auditory system functionality seems paramount. The two-point method discussed above could not be related to the BM dynamic behavior as there is only one source of excitation present, i.e., the incoming sound signal. Therefore, a single-point excitation method is of particular interest in the current work. The research concerning this approach is considerably newer and, in turn, rare compared to the two-point excitation method. In the single-point excitation approach, one of the two point forces would be replaced with an absorber, e.g., a spring-dashpot pair [30-34], that generates the counterforce required to cancel the reflection of the incident wave. In one of the latest published works on this topic, Cheng et al. [34] studied the wave separation phenomenon (i.e., having concurrent TW and SW in separate sections of the structure) in Euler-Bernoulli beams. They showed that the stiffness and damping coefficients could be tuned for a given excitation frequency such that wave separation along the structure would be observed. The non-existence of pure TWs was attributed to evanescent waves contributing to the beam's overall response. Motaharibidgoli et al. [32] experimentally generated concurrent nearly pure TW and SW sections in a beam using a single-point excitation in conjunction with a passive absorber.

Thus, the focus of this work is two-pronged. The BM in the inner ear is capable of exhibiting a non-reflective traveling wave (NRTW) behavior for a wide excitation frequency band. Inspired by this, the first focus is to incorporate some of the characteristics of the biological BM, e.g., varying stiffness properties, to test if it improves the robustness of the model as proposed in [32] to varying excitation frequencies. Therefore, the stiffness of the model is allowed to be distributed along the length of the beam through geometrical tapering or material gradation. The second prong of the present research is to perform a thorough parametric study to understand the model's sensitivity to parameter alteration. These include the absorber's location, excitation frequency, absorber's coefficients, tapering/grading aggressivity, and the gradient index (the order of the power-law formulation used for tapering/grading).

## 2  Problem Formulation and Solution Methodology

The schematic of the general system under study is presented in Figure 1, where $L$ and $L_1$ denote the beam's length and absorber's location, respectively.

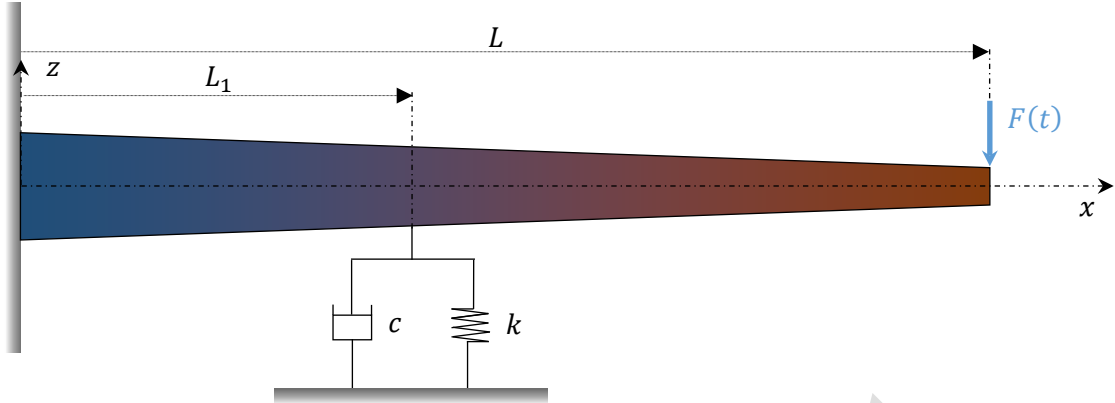

Figure 1. An Euler-Bernoulli cantilever non-uniform non-homogeneous beam with an intermediate passive absorber.

Also, $c$, $k$, and $F$ are, in turn, the absorber's damping and stiffness coefficients and external sinusoidal harmonic force. The general assumption for the model is that any given mechanical or geometrical property variation is mathematically idealized using the power-law formulation [35, 36] given by,

$$P(x) = P_l + \left(\frac{x}{L}\right)^N (P_r - P_l), \qquad (1)$$

where $P(x)$ is a mechanical or geometrical property's value in a lengthwise location. $P_l$ and $P_r$ are associated with the value of said property at the leftmost ($x = 0$) and rightmost ($x = L$) sides of the beam, respectively. $N$ is the gradient index. In this work, the properties that vary according to Eq. (1) are considered to be the thickness $h(x)$, width $b(x)$, Young's modulus $E(x)$, and material density $\rho(x)$.

To develop the mathematical model of the system depicted in Figure 1, Hamilton's principle will be used in the following form [37],

$$\delta \int_{t_1}^{t_2} (T - V + W_{nc}) dt = 0, \qquad (2)$$

where, $T$, $V$, and $W_{nc}$, in turn, denote kinetic and strain energies and the virtual work done by the non-conservative forces. Considering Euler-Bernoulli's beam assumptions, the beam's displacement field is defined as follows,

$$\begin{aligned} U(x,z,t) &= -zw(x,t)_{,x} \\ W(x,z,t) &= w(x,t) \end{aligned}, \quad (3)$$

where, $U$ and $W$ denote the longitudinal and lateral displacements at any given point of the beam in terms of its neutral axis's lateral displacement ($w$). In the notation implemented to develop the formulation in this work, "$,x$" and "$,xx$", in turn, stand for the first and second-order partial differentiation with respect to $x$. Therefore, the beam's strain field is obtained as follows,

$$\varepsilon_{xx} = U_{,x} = -zw_{,xx}. \quad (4)$$

The kinetic and potential energy expressions could be derived in the following general forms,

$$T = \int_D \frac{1}{2}\rho(x)\dot{W}^2 dD = \frac{1}{2}\int_0^L b(x)h(x)\rho(x)\dot{w}^2 dx, \quad (5)$$

and,

$$V = \int_D \frac{1}{2}E(x)\varepsilon_{xx}^2 dD + \int_0^L \frac{1}{2}kw^2\delta(x-L_1)dx = \int_D \frac{1}{2}\left(E(x)\frac{b(x)h^3(x)}{12}w_{,xx}^2 + kw^2\delta(x-L_1)\right)dx. \quad (6)$$

It is worth noting that the rotational kinetic energy is disregarded in Eq. (5) due to Euler-Bernoulli's beam theory assumptions. Finally, the non-conservative virtual work can be calculated as follows,

$$\delta W_{nc} = -\int_0^L \left(c\dot{w}\delta(x-L_1) - F\delta(x-L)\right)\delta w dx. \quad (7)$$

In Eqs. (6) and (7), $\delta(x - L_i)$ denotes the Dirac delta function, which is used to describe point forces in a distributed form in order to enable their incorporation in the beam's dynamic equation. The Dirac delta should not be confused with the variation symbol that for the sake of clarity is implemented here in a non-italicized form (δ). By substituting expressions (5)-(7) in Eq. (2), and with some mathematical manipulation, the equation of motion of the given system is obtained as follows,

$$\rho(x)b(x)h(x)\ddot{w} + \left(E(x)\frac{b(x)h(x)^3}{12}w_{,xx}\right)_{,xx} + kw\delta(x-L_1) + c\dot{w}\delta(x-L_1) = F\delta(x-L). \tag{8}$$

In order to solve Eq. (8) and discretize the problem, the Galerkin method is implemented according to [38]. Therefore, a solution in the following form is assumed,

$$w = \sum_{i=1}^{n} \eta_i(t)\varphi_i(x), \tag{9}$$

where $\eta_i$ is the generalized coordinate and $\varphi_i$ is the trial function. Also, $n$ corresponds to the number of assumed trial functions. For the Galerkin method solution to converge in the present study, a high number of trial functions are required. Therefore, implementing the traditional form of the uniform Euler-Bernoulli beam's mode shape functions may not be feasible as they are known to be inaccurate in the higher modes. For this reason, the modified uniform Euler-Bernoulli eigenfunctions given by Gonçalves et al. [39] are utilized in this study. Inserting the solution given in (9) to the equation of motion (8) and implementing the Galerkin procedure will yield the discretized form of the equation of motion as follows,

$$\mathbf{M}\ddot{\boldsymbol{\eta}} + \mathbf{C}\dot{\boldsymbol{\eta}} + \mathbf{K}\boldsymbol{\eta} = \mathbf{q}, \tag{10}$$

where, **M**, **C**, and **K** denote the mass, damping, and stiffness matrices, respectively. **q** is the forcing vector, and **η** is the generalized coordinates vector. The expressions for the system's matrices and the forcing vector are given as follows,

$$\begin{aligned}
\mathbf{M} &= \int_0^L \rho(x) A(x) \varphi_i \varphi_j \, dx \\
\mathbf{K} &= \int_0^L \left[ \left( E(x) I(x) \varphi_{i,xx} \right)_{,xx} + k \varphi_i \delta(x - L_1) \right] \varphi_j \, dx \\
\mathbf{C} &= \int_0^L \left[ c \varphi_i \delta(x - L_1) \right] \varphi_j \, dx \\
\mathbf{q} &= \int_0^L \left[ F \delta(x - L) \right] \varphi_j \, dx
\end{aligned} \quad , \qquad (11)$$

where, $A(x) = b(x)h(x)$ and $I(x) = b(x)h(x)^3/12$. It is worth noting that these relations are derived in their most generic form, however, special cases where all parameters but one are kept constant will be studied to isolate the effects of that parameter on the overall response.

## 3 Convergence study and code verification

First, a convergence study will be conducted to determine the sufficient number of trial functions for the desired accuracy. Then the code will be verified against previously published results to prove its validity.

### 3.1 Convergence study

As the first step, a uniform homogenous beam with the properties given in Table 1 is considered for the convergence study. The absorber is assumed to be located at $L_1 = 0.4L$.

Table 1. The properties of the considered beam for the convergence study.

| Parameter | Value |
| --- | --- |
| $L$ | $2000\ mm$ |
| $L_1$ | $800\ mm$ |
| $b_l = b_r = b$ | $30\ mm$ |
| $h_l = h_r = h$ | $10\ mm$ |
| $E_l = E_r = E$ | $71\ Gpa$ |
| $\rho_l = \rho_r = \rho$ | $2700\ kg/m^3$ |
| $k$ | $5625000\ N/m$ |
| $c$ | $875\ N.s/m$ |

It is worth mentioning that the authors will focus on the frequency range of 300-3400 Hz (the standard telephony band [40, 41]) since speech sound is mainly concentrated in this frequency range. This was the bandwidth of old analog telephones that ensured that the transmission of sound was intelligible at low costs. Therefore, the excitation frequency is assumed to be the standard telephony band's higher limit. The convergence criterion is automatically satisfied for lower excitation frequencies should the response converge at this frequency.

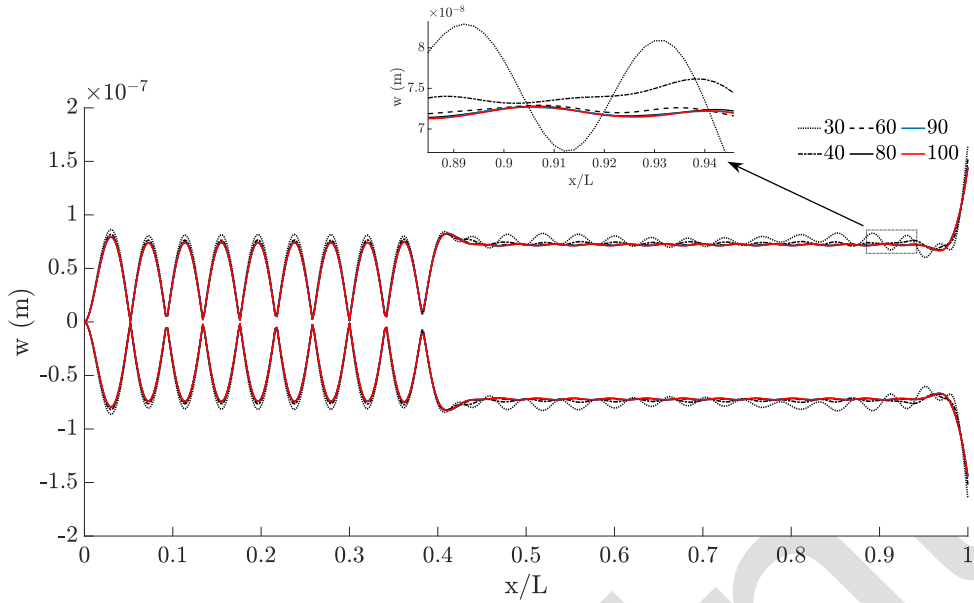

Figure 2. The TW response of a uniform beam using different numbers of trial functions.

The results presented in Figure 2 depict the wave envelopes that are obtained using an increasing number of trial functions, suggesting that the response converges when the number of considered trial functions is increased above 80. Even though using 80 trial functions suffices, the remainder of this study will use 100 trial functions as a prudent approach to compensate for the natural frequency shifts due to parametric alterations such as geometric tapering and material grading. As an example, the case where the beam is aggressively tapered is considered for a convergence study to test whether 100 trial functions would satisfy the convergence criterion in such cases. For this purpose, $L$, $L_1$, $E$, and $\rho$ are kept the same as the previous case, however, it is assumed that the $b_l = 45\ mm$, $b_r = 15\ mm$, $h_l = 15\ mm$, and $h_r = 5\ mm$, implying that the leftmost cross-sectional area of the beam is nine times the right one. The intermediate absorber coefficients are $k = 1875000\ N/m$ and $c = 2187.5\ N.s/m$. The result of this convergence study is illustrated in Figure 3.

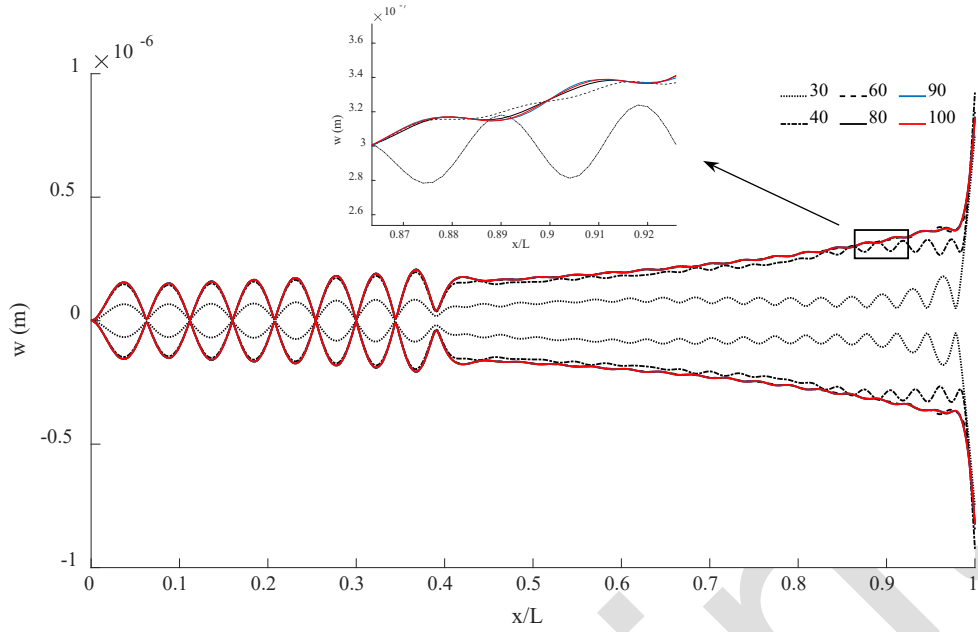

Figure 3. The TW response of the tapered beam using different numbers of trial functions.

Figure 3 suggests that the results are almost perfectly converged at 80 trial functions, meaning that the implementation of 100 trial functions is more than enough to take into account the natural frequency shifts due to parametric changes.

### 3.2 Verification of the developed code

A code based on the presented formulation is developed and verified against some of the results reported in the literature. In the first step, per the study published by Craver and Jampala [42], let $c = 0$ and $L_1 = 0.6L$. The material properties are considered to be homogeneous, i.e., $E_l = E_r = E$ and $\rho_l = \rho_r = \rho$. The taper ratio for thickness and width is 1.4, meaning that each of these dimensions is 1.4 times its value at the fixed end of the beam. The fundamental frequency obtained from the formulation developed herein is validated against the analytical results given by [42] and shown in Table 2.

Table 2. Verification of the first nondimensionalized natural frequency of a tapered beam with an intermediate absorber using [42]

| $kL^3/EI(L)$ | $\omega_1\sqrt{(\rho A(L)L^4)/EI(L)}$ | | Error (%) |
|---|---|---|---|
| | [42] | Present study | |
| 0 | 5.6483 | 5.6478 | 0.01 |
| 1 | 5.7172 | 5.7072 | 0.17 |
| 10 | 6.2943 | 6.2099 | 1.34 |
| 50 | 8.2427 | 7.9580 | 3.45 |
| 100 | 9.9085 | 9.5057 | 4.07 |
| 500 | 14.8324 | 14.4252 | 2.75 |
| 505.8 | 14.8629 | 14.4578 | 2.73 |
| 1000 | 16.3405 | 16.0735 | 1.63 |

As suggested by Table 2, the results are in fair agreement for different values of the absorber's stiffness coefficient. The slight discrepancy arises due to the fact that a semi-analytical method is being compared against the exact solution, thus it can be deemed acceptable.

As the next step, the TW response (the wave envelope) of a uniform homogenous beam obtained by this code will be compared to that of Motaharibidgoli et al. [32]. In this case, $L = 1.5812\ m$, $L_1 = 0.4L$, $E_l = E_r = E = 71\ GPa$, $\rho_l = \rho_r = \rho = 2700\ kg/m^3$, $h_l = h_r = h = 1.5875\ mm$, and $b_l = b_r = b = 12.7\ mm$. The absorber's coefficients are $k = 71\ kN/m$ and $c = 9\ Ns/m$. The beam is excited at $f = 1300\ Hz$. Also, it must be noted that in all the results presented in this study, the forcing amplitude is unity and directed downwards.

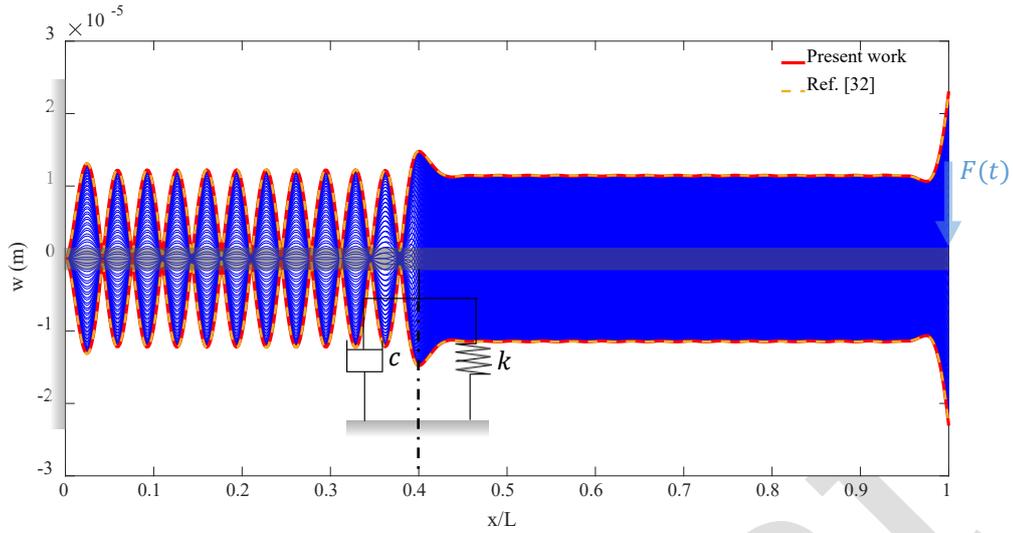

Figure 4. TW response comparison between the present study and [32]

Figure 4 shows perfect agreement between the wave envelopes obtained by the current method and that of [32]. Therefore, the developed code is successfully verified.

## 4 Results and Discussion

The effects of different parameters, i.e., excitation frequency, absorber's location, cross-sectional tapering, material properties gradation ($E$ and $\rho$), and the gradient index on the TW response are studied in the remainder of this paper. As the first step to these investigations, one must establish a method to determine the presence (or the quality) of the TWs in the dynamic response of the beam. The cost function (CF) discussed in [32] is adopted in response to this requirement. It is worth mentioning that in the case of tapered beams, the overall wave envelope will be slightly sloped along the length of the beam due to the lengthwise varying stiffness or mass distribution. Therefore, the mean CF value is considered in such cases, i.e., a mean of the CF calculated using consecutive peaks and valleys along the wave envelope associated with the TW section.

For instance, the CF plot for the system whose TW response is presented in Figure 4 is depicted in Figure 5. Based on the spring-damper system coefficients, a CF value is calculated using the obtained wave envelope with the formulation mentioned above. The value of the CF is somewhere between 0 and 1, where 0 is the indicator of pure TW, and in contrast, 1 corresponds to pure SW. Any value between these numbers can be attributed to the presence of a mixture of both types of the discussed waves in the response, i.e., a hybrid wave. As the CF approaches 0, TWs become more dominant than SWs in the overall response. In the remainder of this study, it is assumed that CF values below 0.3 are desirable, meaning that the TWs mostly dominate the overall response over SWs.

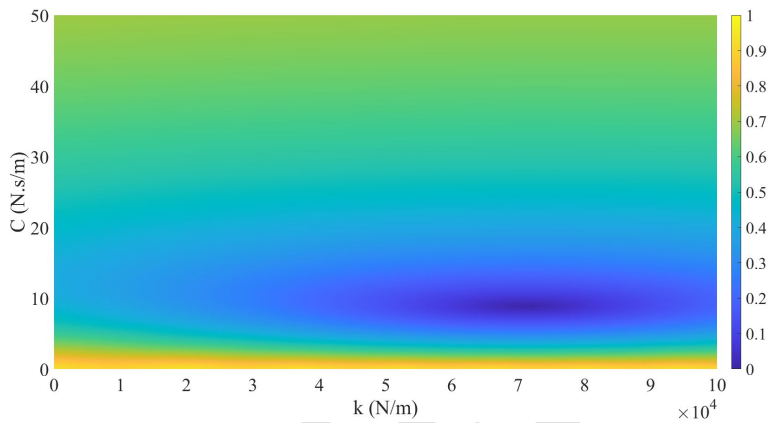

Figure 5. The CF plot for uniform homogeneous beam studied in [32]

## 4.1 Effect of excitation frequency ($f$)

The TW responses of the uniform homogeneous beam as well as the non-uniform and non-homogenous cases within the excitation frequency bandwidth are studied in the following subsections.

### 4.1.1 Uniform homogenous beam

One of the critical observations in the BM's dynamic behavior is the presence of NRTWs in a wide range of excitation frequencies (i.e., human auditory range, roughly 20 Hz-20 kHz [43]). Therefore, the goal is to investigate the robustness of the model presented here to the variation of excitation frequencies within the standard telephony band. First, the plot of stacked CFs associated with the TW response of the uniform beam described in Table 1 for the excitation frequency range is presented in Figure 6(a). The CF plots for each frequency line ($\Delta f = 10\ Hz$) are piled on top of each other to enable one to inspect if a common value (or a region of values) for the absorber's stiffness and damping coefficients exists where the CF has a desirable value (i.e., $CF \leq 0.3$) at all excitation frequencies.

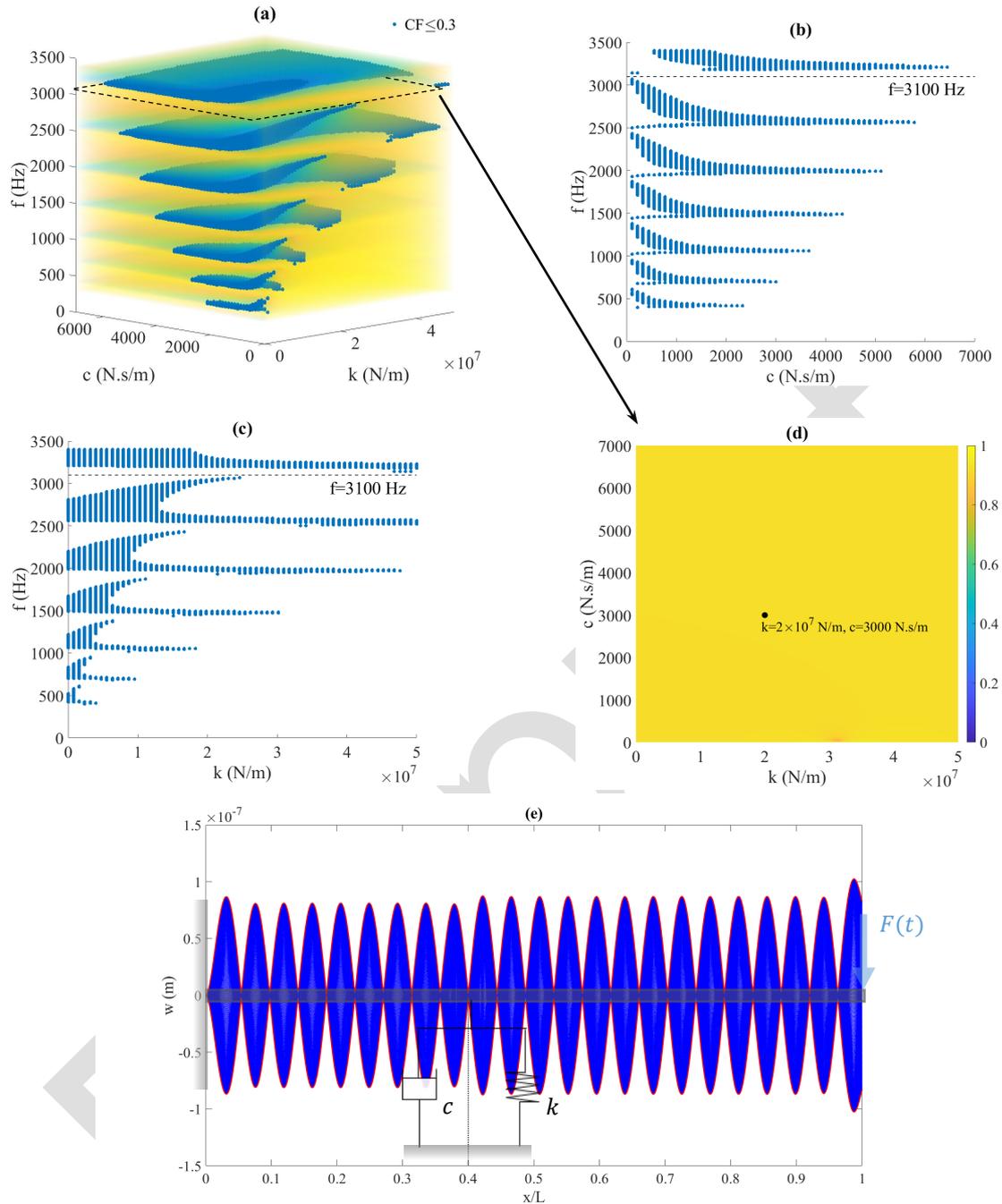

Figure 6. The stacked CF plots for the uniform homogeneous beam with marked below 0.3 values for the standard telephony range in (a) 3D view (b) $f-c$ view (c) $f-k$ view. (d) The CF plot for $f = 3100\ Hz$ (e) The wave response of the beam with the absorber coefficients $k = 2 \times 10^7\ N/m$ and $c = 3000\ N.s/m$.

Figures 6(a) to 6(c) suggest that in the configuration of a uniform homogeneous beam with an absorber at $L_1 = 0.4L$, it is not possible to obtain an NRTW response for all excitation frequencies in the selected band. As one may readily notice, in some frequency lines (slices of the stack plot), no low-valued CF seems to exist. For instance, Figure 6(d) shows the CF plot at $f = 3100\ Hz$ where no optimal value exists at all. In fact, the CF is close to 1 (pure SW) for all $k$ and $c$ values. This suggests that the absorber's location coincides with or closely approaches the vibrational node at these frequencies for all values of the absorber's coefficients in their assumed ranges. Figure 6(e) is the physical representation of this discussion, where for given absorber coefficients picked from the CF plot presented in Figure 6(d), almost pure SWs are observed in both beam regions. On the other hand, the general trend in Figure 6(a) suggests that the desirable CF values region (which resembles a repetitive curved cone pattern) grows larger as the excitation frequency increases, meaning that the CF is desensitized to the variation of absorber's coefficient values. At the next step, a similar investigation is carried out to see if manipulations such as tapering and functionally grading the beams could affect the trend observed in Figure 6.

### 4.1.2 Tapered and graded beams

In Figure 7(a), it is assumed that $b_l = 48\ mm$, $h_l = 16\ mm$, $b_r = 12\ mm$ and $h_r = 4\ mm$, meaning $A(0) = 16A(L)$. In addition, the beam is considered to be linearly tapered, i.e., $N = 1$. Other parameters are kept the same as the uniform homogeneous beam studied earlier. Figure 7(a) suggests that tapering the beam does not lead to the formation of a continuous region along the frequency axis where common optimal values of $c$ and $k$ exist independent of excitation frequency. This is due to the same reason that was observed in the uniform homogeneous beam. Therefore, tapering does not address this issue as suspected.

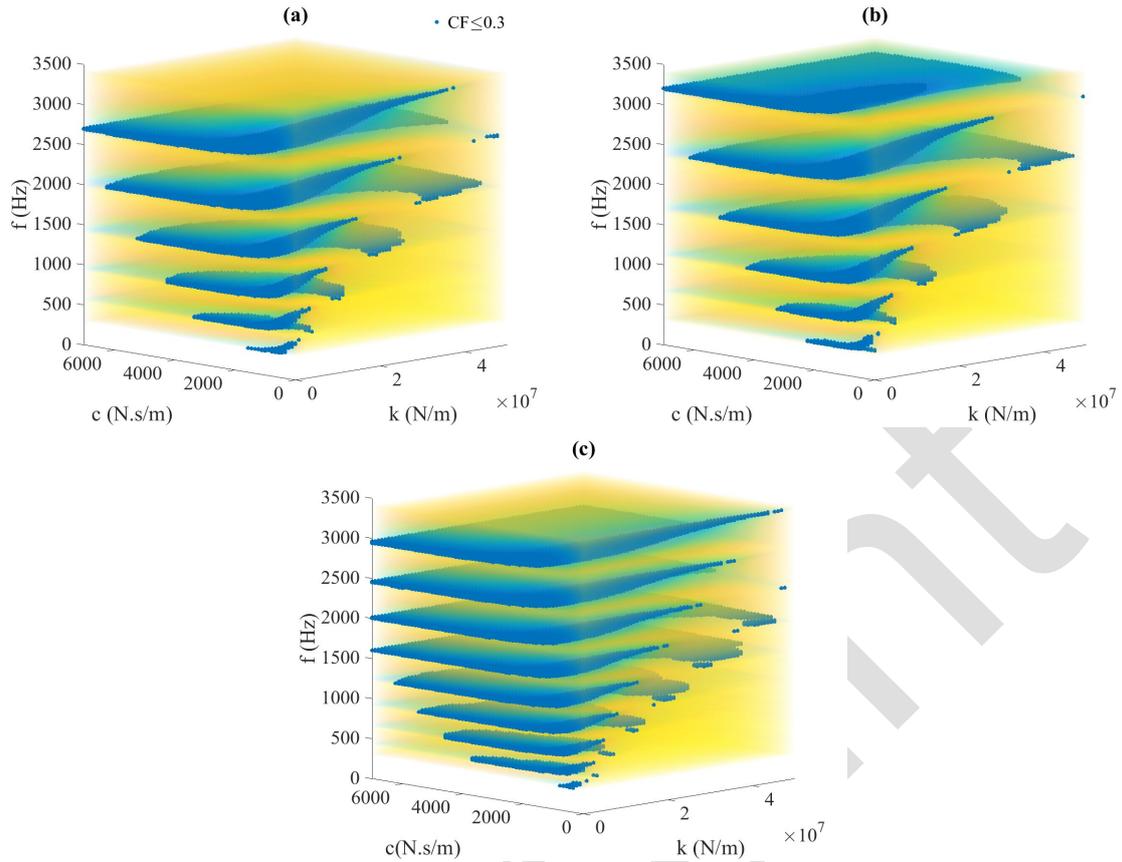

Figure 7. The stacked CF plots for the standard telephony range for (a) the tapered beam, (b) The beam where $E$ is graded, and (c) The beam where $\rho$ is graded.

Figures 7(b) and 7(c), in turn, correspond to the cases where $E$ and $\rho$ are linearly graded. In Figure 7(b) $E_l = 227.2\ GPa$ and $E_r = 14.2\ GPa$, i.e., $E_l = 16E_r$. Also, in Figure 7(c) $\rho_l = 8640\ kg/m^3$ and $\rho_r = 540\ kg/m^3$, i.e., $\rho_l = 16\rho_r$. It is worth mentioning that these ratios are chosen as examples of aggressive tapering/grading cases. In both cases, it is assumed that the remaining parameters are the same as the uniform homogenous case. Again, similar trends to that of the uniform homogeneous beam are observed in both cases. However, grading the density seems to increase the frequency-wise instances where optimal values exist. This must be attributed to the general decrease in the natural frequencies of the system due to the beam's density grading, which gives rise to the contribution of higher modes in the beam's response within the assumed excitation frequency range. In this case modes 9 to 27 of the graded beam lie within the excitation band, whereas in the uniform homogeneous case, modes 8 to 24 are within the range. This means that as we sweep through the excitation frequency band, the distance between vibration nodes is closer in the graded beam case, thus the coincidence of nodes and absorber's location occurs more frequently. It can be concluded that since the coincidence of the absorber's point of attachment with the nodal points is inevitable when continuously sweeping through the excitation frequency band, it may not be possible to obtain a continuous 3D manifold of optimal CF values using this configuration.

## 4.2 Effects of the absorber's location

Here, it is assumed that the absorber's location is gradually shifting within the range of $L_1 = [0.1L, 0.75L]$. The beam is excited at $f = 3400\ Hz$, and all other material and geometric properties are in accordance with the uniform homogeneous beam discussed in the previous section.

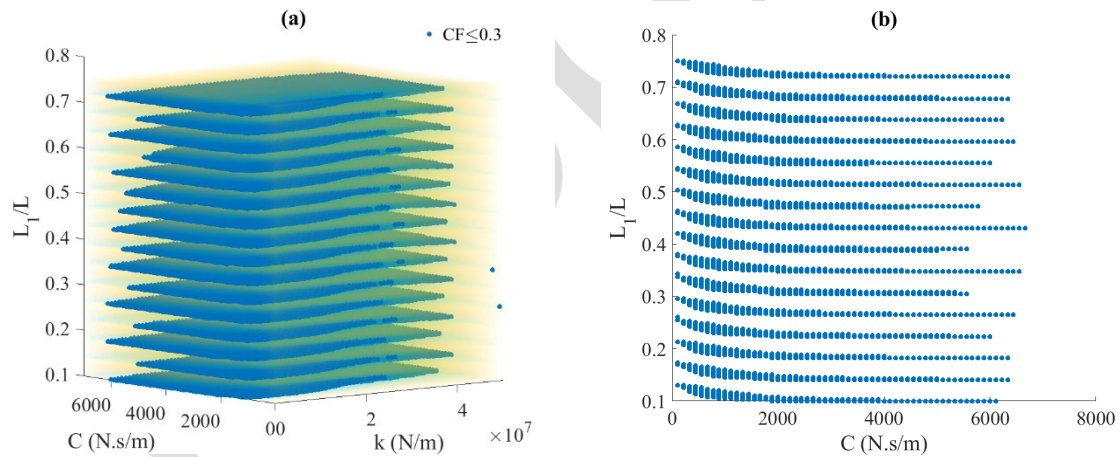

Figure 8. The CF plots stack for the uniform homogenous beam where absorber's location changes. (a) 3D view and (b) side view.

The trend observed in Figure 8(a) makes sense since as the absorber's location moves along the length of the beam it coincides with the vibration nodes. The absorber is expected to become neutralized in these locations, and only a SW response will be observed in the beam. As the absorber's location approaches the node, the optimal CF region shrinks and gradually disappears, and when that node is passed, it starts to expand in size again. Therefore, when a TW response is desired based on an excitation frequency, one should avoid the corresponding vibration nodes when picking a location for the absorber. Furthermore, Figure 8(b) suggests each optimal region cone located between two vibration nodes tends to become narrower as the absorber moves toward the following node in $+x$ direction. No matter which point is chosen for installing the absorber, it is going to be between two vibrational nodes, thus, it is suggested to pick this location in a manner that it is closer to the left-hand side node where the base of each optimal CF cone is located. This way, there would be more flexibility in terms of choosing an absorber.

### 4.3 Effects of cross-sectional tapering on the beam's TW response

Next, the effect of tapering on the beam's response in a given excitation frequency when all other parameters are kept constant is studied.

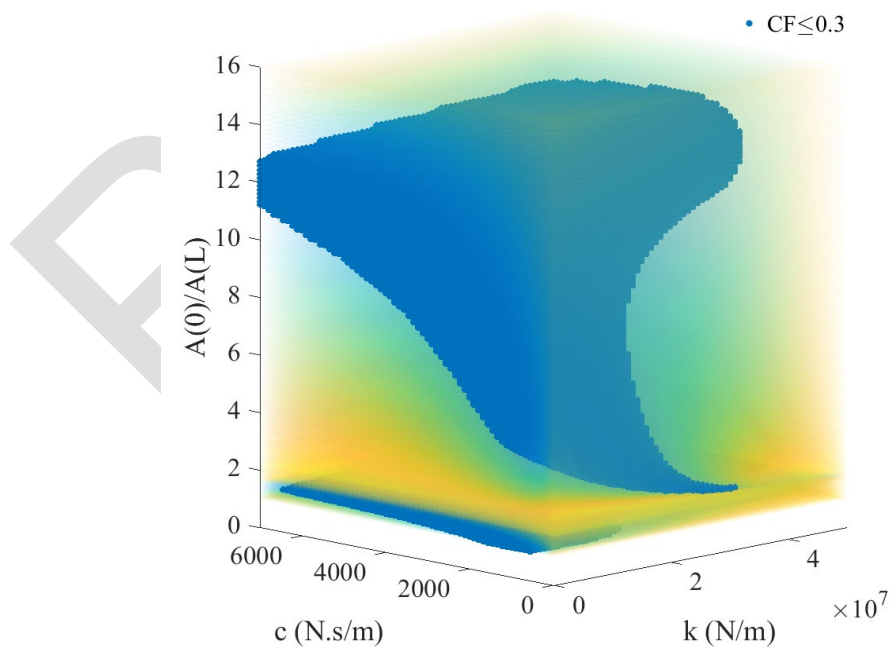

Figure 9. The CF plots stack for varying taper ratios.

In Figure 9, it is assumed that the structure is excited at $f = 3400\ Hz$. The beam will be gradually tapered from uniformity to the point where $b_l = 48\ mm$, $h_l = 16\ mm$, $b_r = 12\ mm$ and $h_r = 4\ mm$, i.e., $A(0) = 16A(L)$. This results in an increase in the optimal CF region area per taper ratio, which should be attributed to the vibrational nodes' movement relative to the absorber's location due to gradual natural frequency shifts and slight mode changes as a result of the system's stiffness change. These node shifts are shown in Figure 10, where the responses of the beams with area ratios of 6 and 12 are given, respectively. The absorber coefficients are $k = 7500000\ N/m$ and $c = 2400\ N.s/m$ and all other parameters are given in Table 1. Also, $f = 3400\ Hz$. It is readily observable that the depicted node gets closer to the absorber's location as the ratio increases.

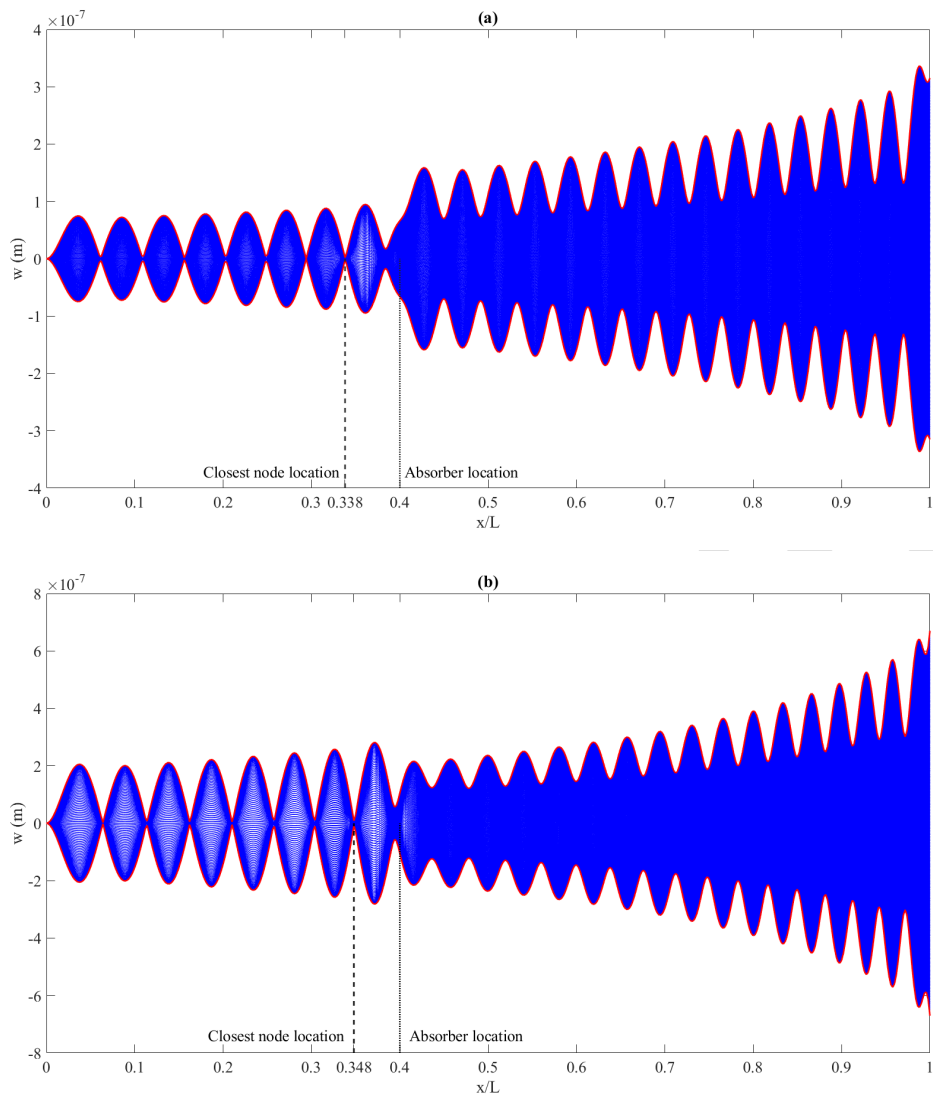

Figure 10. The displacement response of cross-sectionally tapered beams with (a) $A(0)/A(L) = 6$ and (b) $A(0)/A(L) = 12$ ratios

Suppose the ratio continues to increase beyond 16. In that case, the pattern observed in Figure 9 is expected to repeat itself as the absorber's location falls within the following pairs of nodes. Therefore, between any two nodal points, as the ratio increases, the optimal region area grows larger. It is worth reemphasizing that here the spring-damper system is fixed, and the nodal locations are shifting due to varying taper ratios.

## 4.4 Effects of grading material properties on the beam's TW response

In this section, the effects of grading $E$ and $\rho$ throughout the beam's length are investigated. Again, it is assumed that in each case, all other material properties are identical to those of the uniform homogeneous case. The beam is again excited in $f = 3400\ Hz$, and the material properties ratio of interest would vary from 1 to 16 for both cases.

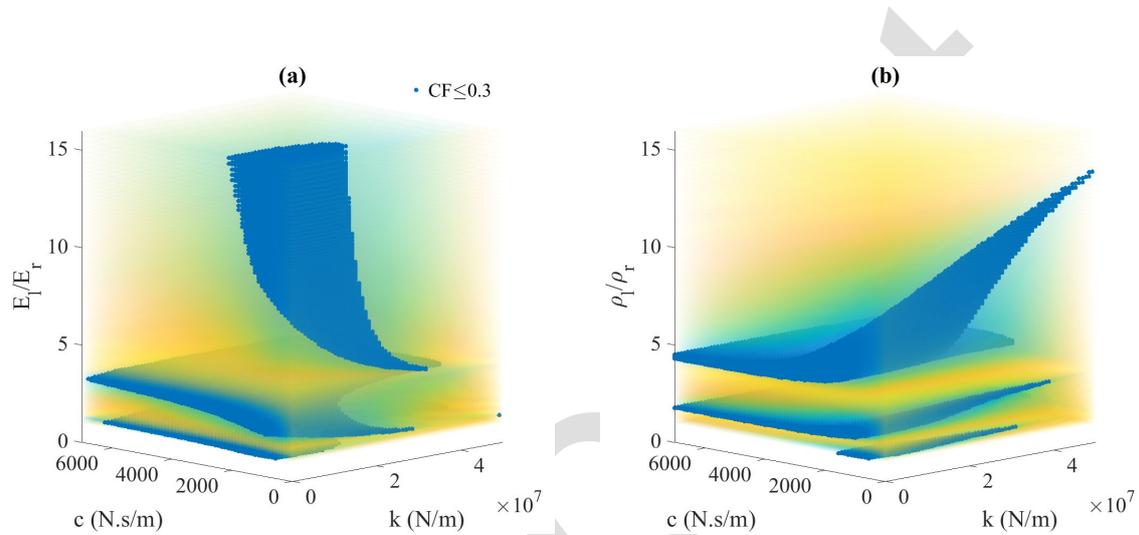

Figure 11. The CF plots stack for: (a) Varying $E$ ratio and (b) Varying $\rho$ ratio cases

Figure 11(a) reveals that the effect of grading the modulus of elasticity is similar to tapering as they have a similar effect on the beam's stiffness. On the other hand, a reversed trend is observed when the grading is performed on the structure's material density, shown in Figure 11(b). This makes sense as density is directly associated with the structural mass rather than stiffness; therefore, such a reverse effect is expected. Therefore, similar to previous cases, gradual grading and tapering would lead to natural frequency shifts, giving rise to the occasional coincidence of vibration nodes with the absorber's location. This would lead to the repetitive patterns observed in the presented results. Moreover, it can be seen that the rate of change is significantly higher in smaller density or modulus of elasticity ratios. It should be noted that each of these cones is associated with the relative location of the absorber within two adjacent vibrational nodes. Therefore, the results presented in Figure 11 suggest that in lower grading ratios the absorber's relative location changes from one node to the next one at a faster rate than higher ratios. For instance, Figure 12 shows the relative nodal locations with respect to absorber location and will be used to elaborate on the phenomenon mentioned above.

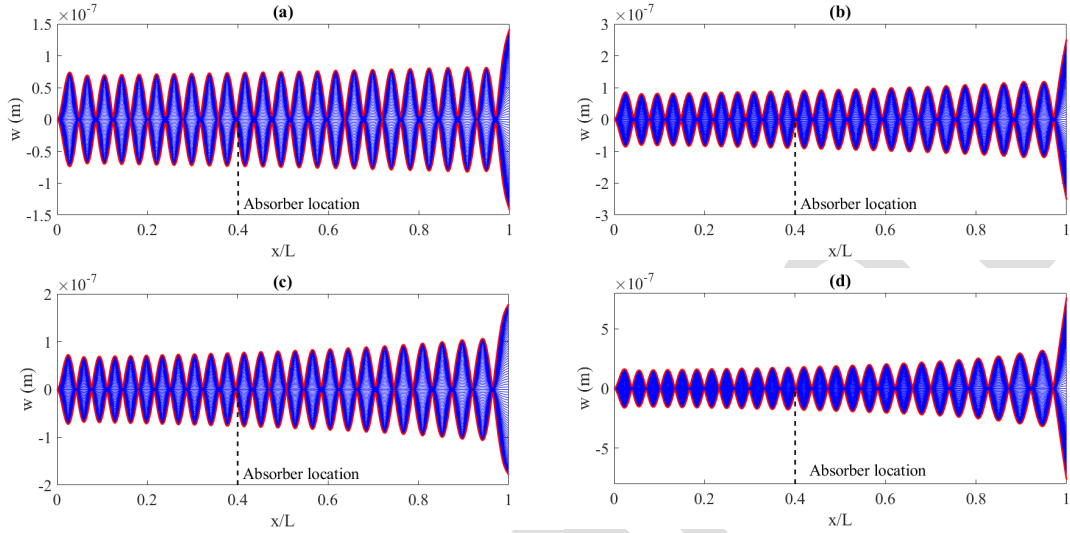

Figure 12. The nodal locations relative to the absorber location for cases where (a) $\rho_l/\rho_r = 1.7143$, (b) $\rho_l/\rho_r = 3.5674$, (c) $\rho_l/\rho_r = 4.1637$, and (d) $\rho_l/\rho_r = 13.8723$

Figure 12 shows the response of the beam where the density is graded with different ratios, i.e., 1.7143, 3.5674, 4.1636, and 13.8723, respectively. The absorber is turned off for these cases to solely inspect the effect of grading in vibrational node shifts. All other parameters are set according to Table 1 and $f = 3400\ Hz$. Figure 12(a) depicts that node 11 is coinciding with $x = 0.4L$, whereas in Figure 12(b) where the grading has slightly become more aggressive, node 12 is located at about $x = 0.4L$. The same shift is shown from node 12 to 13 in Figures 12(c) and 12(d), respectively for a considerably bigger change in grading ratio. This suggests that the same absorber location shift within the range of two consecutive pairs of nodes is happening at different rates of grading ratio changes. The first transition takes place as $\rho_l/\rho_r = 1.7143$ is increased to $\rho_l/\rho_r = 3.5674$, which is shown by the shorter cone in Figure 11(b). The second one that is associated with the taller cone is slower and takes place as the ratio increases from 4.1637 to 13.8723.

## 4.5 Effects of altering the gradient index (N)

Up to this point, only linear tapering and gradation are considered, meaning that in Eq. (1), $N = 1$. Here, the goal is to investigate the effect of using values other than 1 for $N$, in the TW response of the tapered/graded beam. This provides a more aggressive change of properties in different portions of the beam with length depending on $N$'s value.

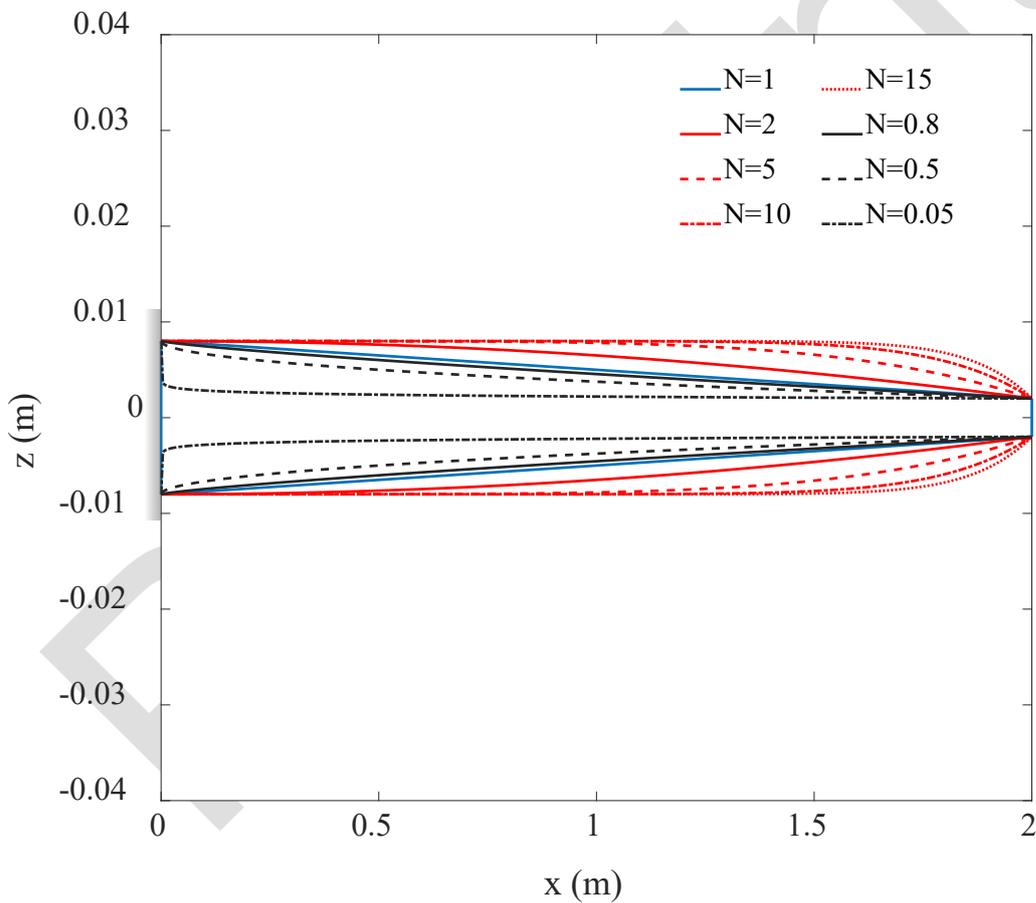

Figure 13. The side view of the tapered beam for different values of gradient index ($N$)

Figure 13 provides the side view of the cantilever tapered beam for different values of $N$. As one could see the rate of change is higher for values below $N = 1$. This suggests that the change in the optimal CF region area should be more gradual for values above 1. This is confirmed by the results provided in Figure 14.

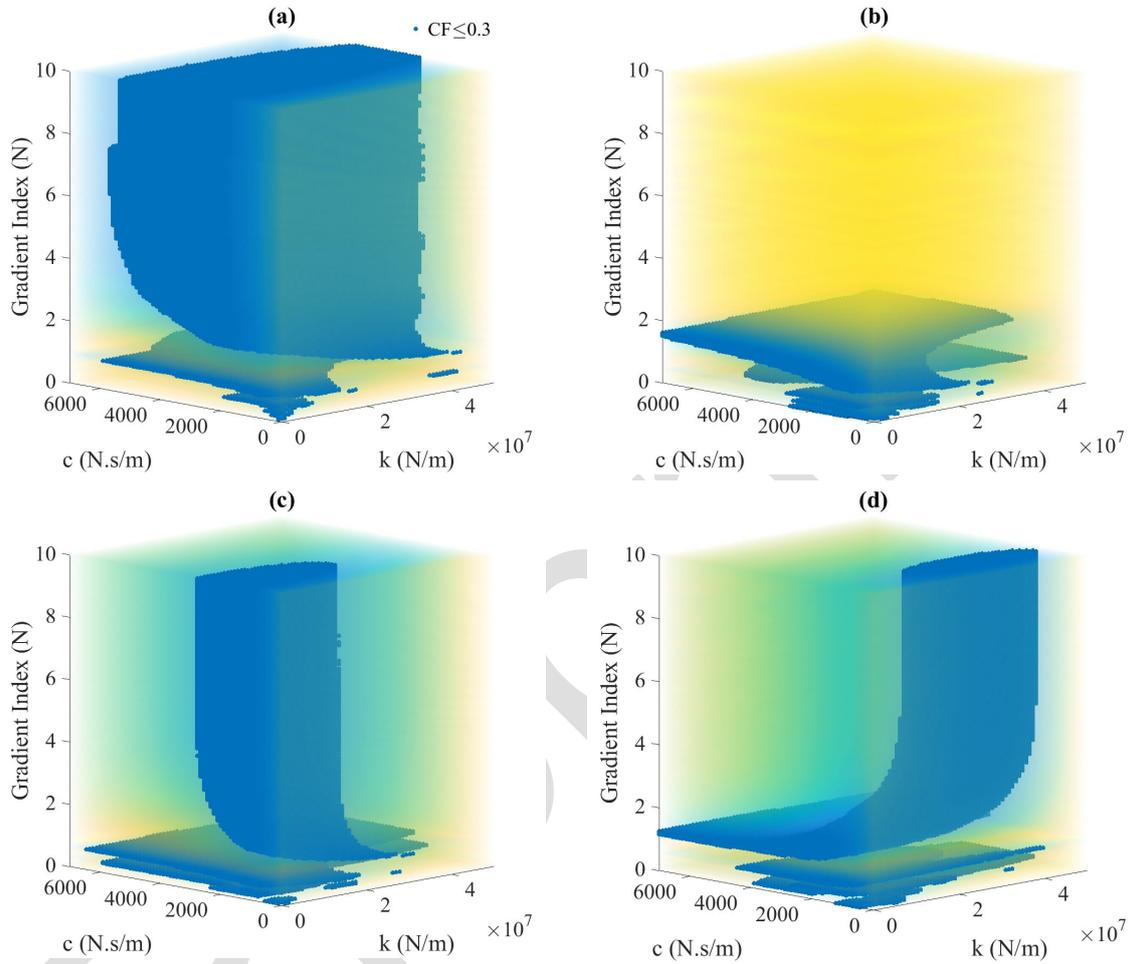

Figure 14. CF plots stack for varying $N$ for (a) the tapered beam, (b) graded $E$ and $L_1 = 0.4L$, (c) graded $E$ and $L_1 = 0.42L$, and (d) graded $\rho$.

Figure 14(a) is associated with the homogeneous tapered beam where, $b_l = 48\ mm$, $h_l = 16\ mm$, $b_r = 12\ mm$ and $h_r = 4\ mm$. The absorber is located at $L_1 = 0.4L$ and it is being excited at $f = 3400\ Hz$, similar to the previous cases. As one may readily notice, increasing the gradient index in the case of the tapered beam can lead to an increase in the robustness of the model in terms of acceptable absorber coefficients values. However, this must be attributed to the fact that this is due to its effect on stiffness properties and the consequent nodes' shifts relative to the absorber's location. At roughly $N = 4$, a point is reached wherein the CF optimal region area would no longer change in response to increasing the gradient index. This is of particular interest in the case of manufacturing since a more aggressive change would create manufacturing complications for little to no return in performance. A similar trend is also observed in cases where the structure undergoes material gradation. This saturation of sorts is more notable in Figure 14(b), where $E_l = 227.2\ GPa$ and $E_r = 14.2\ GPa$, and the absorber and excitation circumstances are the same as in the previous case. In this case, as $N$ gradually increases a vibration node seems to coincide with the absorber's location. Therefore, at $N \cong 2$ the existence of an optimal CF region disappears. As $N$ further increases, it no longer affects the beam's modal properties, and therefore, the coinciding node would not move from the point where the absorber is located. As shown in Figure 14(c), with similar conditions and moving the absorber to $L_1 = 0.42L$, a similar behavior to tapering effect could be observed. Figure 14(d) corresponds to the case where the beam's density is graded with $\rho_l = 8640\ kg/m^3$ and $\rho_r = 540\ kg/m^3$. The saturation is again observed in this case however the trend is expectedly the reverse to that of tapering and elasticity modulus gradation, in the sense that the optimal region firstly shrinks as $N$ is increased and then remains constant as it further increases. Therefore, in general, manipulating $N$ could be an effective approach in broadening the optimal region area. However, after the discussed saturation takes place, no significant change is obtained in the optimal CF region area by further increasing $N$.

## 5 Conclusion

Inspired by the vibratory trend observed in the mammalian BM, the TW formation in a thin generic Euler-Bernoulli beam coupled to a passive spring-damper system as the reflected wave absorber was studied. As the first step, based on the Euler-Bernoulli beam assumptions, using Hamilton's principle, the lateral displacement dynamic equation of motion in its most general form was developed. Then the Galerkin method was applied to the equation to develop the discretized form of the problem. The uniform homogeneous beam's modified eigenfunctions were used as trial functions fed to the Galerkin algorithm. A code was developed using the derived formulation and was successfully verified against the published works available in the literature. After the required number of trial functions for the solution to properly converge within the frequency range of interest was determined, a full-scale parametric study was conducted to investigate the effect of different system parameters on the beam dynamic response. The concluding remarks of this study can be listed as follows:

- Based on the results obtained in this study, it proves challenging to obtain a high-quality NRTW response of the beam for a relatively wide range of excitation frequencies with a combination of a beam coupled to a single spring-damper system. This is because of the geometrical limitations that are imposed by the nodal points of vibration.
- The gradual alteration of any parameter, e.g., the beam's stiffness through tapering or elasticity modulus grading, would affect the system's natural frequencies, which, in turn, gives rise to the vibration nodes being shifted along the beam continuously for a given excitation frequency. Therefore, in specific instances, the nodes will coincide with the absorber's location, neutralizing its absorbing effect. As a result, in these instances, only SWs will be observed on the beam. In this simplified model of the BM, these instances always exist for any given combination of these structural elements.
- Even though a continuous 3D manifold of optimal CF values seems nonexistent in this case, local improvements were observed when tapering and elasticity modulus grading ratios were increased between two nodal points coincidences.

- The same improvement can be achieved by decreasing the density ratio in a similar situation.
- The Gradient index is observed to be asymptotically ineffective after increased up to a certain point in each case, be it tapered or graded. Therefore, in terms of effectiveness, an optimal $N$ exists for each case, where increasing the order of Eq. (1) beyond such value does not affect the response any further and would highly complicate physical realization through manufacturing.

One of the main motivations of the present study was to investigate the model's capability to passively cancel the incident wave reflection for a wide excitation frequency range. The present study suggests that the current model may not serve as a feasible solution to this problem, as the existence of the nodes is a geometrical limitation that is not easy to overcome with a single spring-damper system. Therefore, as one of the next steps in this work, the authors aim to replace the absorber system presented herein with other types of passive absorbers with multi-mode capabilities to develop a robust TW model.

## Acknowledgments

The authors would like to acknowledge Texas A&M High Performance Research Computing (HPRC) for providing the advanced computing resources needed for this study. The authors would also like to recognize the support provided by Texas A&M J. Mike Walker '66 Department of Mechanical Engineering through Byron Anderson '54 Fellowship and Graduate Summer Research Grant.